%% file: connections-of-dh.tex
\documentclass[a4paper,10pt]{amsart}
\usepackage[utf8]{inputenc}
\usepackage[T1]{fontenc}

\usepackage{graphicx}

\usepackage{amsmath}
\usepackage{amssymb}
\usepackage{amscd}
\usepackage{indentfirst}

\usepackage[pdfpagelabels]{hyperref}

\input{math-lib.tex}

\usepackage{dirtytalk}

\newcommand{\dbar}{\overbar{\partial}}
\newcommand{\polynomials}{\operatorname{Poly}}
\newcommand{\projectivespeciallinear}{\operatorname{PSL}}
\newcommand{\speciallinear}{\operatorname{SL}}
\newcommand{\tracesquared}{\operatorname{tr}^2}

\begin{document}

\title{Some connections of complex dynamics}
\author{Alexandre De Zotti}

\begin{abstract} \ 
 We survey some of the connections
 linking complex dynamics to other fields of mathematics and science.
 We hope to show that complex dynamics is not just interesting on its own but also has value
 as an applicable theory.
\end{abstract}

\maketitle

\section{Introduction}

\emph{Complex dynamics} is the study of the iterations of holomorphic maps
\footnote{There is also another unrelated field called \textit{Complex dynamics} which can also be described as ``nonlinear dynamics''
and usually involves the coupling of different systems, hence the use of the word ``complex''.
In our perspective the word ``complex'' is to be understood as relating to the complex numbers.}
and the field of \emph{dynamics in one complex variable}
is the subfield of complex dynamics concerning the iteration of holomorphic functions defined on Riemann surfaces.
Usually this involves an open connected subset $ U $ of a Riemann surface such as the complex plane $ \complexnumbers $
or the Riemann sphere $ \riemannsphere $
and a non constant holomorphic function $ f $ defined on  $ U $ and whose range intersects $ U $.
Hence if one has to explain to non experts what complex dynamics is then one would have to explain many concepts and ideas:
complex numbers, holomorphic functions, iteration and finally, why studying complex dynamics.

The study of the iteration of holomorphic function belongs to both the fields of complex analysis and the theory of dynamical systems.
 Because of the rigidity of holomorphic functions, the theory of complex dynamics is rich in deep results:
 a complete combinatorial description of the  structure of the Julia set of polynomials and (conjecturally) of the Mandelbrot set
 \cite{Milnor2000},
 application of the
 thermodynamical formalism to Julia sets
 (see for example the survey \cite{Przytycki2018ax})
 which allows to compute their Hausdorff dimensions,
 interplay with circle map dynamics and the theory of small divisors
 (for example \cite{Yoccoz2002}),
 as examples of realizations of unusual topologies as Julia sets
 (e.g. \cite{BlokhBuffCheritatOversteegen2010,RempeGillen2016ax})
 or of pathological dynamical systems \cite{BuffCheritat2012}.

In this article we will attempt to give some examples of what makes complex dynamics an attractive field of research in the point of view of applications.
We do not claim to be exhaustive.
Section \ref{sec: background material} contains some background material on complex dynamics.
In the following sections we will focus on three areas: Kleinian groups, root finding algorithms and the Ising model.
 For each area we will try to explain some of their relations with the field of complex dynamics
 and will refer to further references for more in depth exploration.
Finally in Section \ref{sec: other} we will give quick indications about other connections.

The author wishes to thank Prof. Kuntal Banerjee for helpful discussions.

\section{Background material in complex dynamics} \label{sec: background material}

In this section we cover basic material about complex dynamics.
For some references on the topic see for example \cite{Milnor2006book}, \cite{McMullen1994book}, \cite{CarlesonGamelin1993} or \cite{Beardon1991}.

The most elementary type of holomorphic functions to study in complex dynamics would be polynomials.
 But beyond the trivial case of affine maps, the theory of polynomial dynamics is already rich and complex.
 Let $ d \geq 2 $ and denote by 
 \begin{equation}
  \polynomials (d) 
  \end{equation}
  the set of polynomials of degree $ d $.
 Let $ P \in \polynomials (d) $.
 Then for each $ z = z_0 \in \complexnumbers $
 we define inductively its \emph{orbit} $ \sequence{ z_n }_n $ under $ P $ as
 \begin{equation}
  z_{ n + 1 } \defeq P (z_n) = P^n (z_0),
 \end{equation}
 where $ P^n = P \composedwith \dots \composedwith P $ denotes the \emph{$ n $th iterate} of $ P $.

 Using a direct computation, it is easy to show that if $ \abs{ z } $ is large enough then its orbit
  converges quickly to $ \infty $. This motivates the definition of the \emph{filled Julia set} $ K (P) $ of $ P $:
  \begin{equation}
   K (P) \defeq \setof{ z \in \complexnumbers : \mbox{ the orbit of } z \mbox{ is bounded } }.
  \end{equation}
  Then the \emph{Julia set} $ \juliaset (P) $ is defined as the boundary of the filled Julia set.
  The sets $ K (P) $ and $ \juliaset (P) $ are both totally invariant, that is
  $ P (K(P)) = P^{ -1 } (K(P)) $ and $ P (\juliaset (P)) = P^{ -1 } (\juliaset (P)) $.
  If follows from its definition that the Julia set is characterized by sensible dependence on initial conditions on its neighborhood.
 We will later see a more general definition for the Julia set.

  For example the Julia set of the map $ z \mapsto z^2 $ is simply the unit circle $ \setof{ z: \abs{ z } = 1 } $
  and its filled Julia set the closed unit disk.
  Orbits inside the open disk are attracted by the fixed point at $ 0 $
  and orbits outside the closed unit disk diverge quickly to $ \infty $.
  The unit circle is situated at the interface between these two very distinct behaviors.
 
  The Julia set of a polynomial is either connected or consists of uncountably many connected components.
  And in the latter case
 when $ d = 2 $,
 the Julia set is a Cantor set.
 
 A \emph{critical point} of $ P $ is a point where the derivative of $ P $ vanishes.
  One of the main principles of complex dynamics is that the orbits of the critical points
  determine the global features of the dynamics of the map.
 This principle is exemplified in the following equivalence:
 the Julia set of a polynomial is connected if and only all the critical points of $ P $
 belong to the filled Julia set.

It is easy to see that any polynomial of degree $ 2 $
 is conjugated via an affine change of variables to a polynomial of the form
 \begin{equation}
  P_c (z) = z^2 + c
 \end{equation}
 where $ c \in \complexnumbers $ is some complex parameter.
 Moreover if $ c \neq c' $ then $ P_c $ and $ P_{ c' } $ are not affinely conjugated.
 The family of polynomials $ P_c $ is called the \emph{quadratic family}.
 The quadratic family encompasses the dynamics of all the quadratic polynomials up to affine change of variables.
 The set of parameters $ c \in \complexnumbers $ is also called the \emph{parameter space} of the quadratic family.
 This is just the complex plane $ \complexnumbers $ seen as a family of distinct dynamical systems.
 
In general one studies the properties of a family of holomorphic maps, such as bifurcations, inside the parameter spaces.
 For example the \emph{connectedness locus} of a parametrized family is the set of parameter for which the Julia set is connected.
 The connectedness locus of the quadratic family is more famously known as the \emph{Mandelbrot set}.

As mentioned in the introduction the theory of complex dynamics is concerned with any type of holomorphic functions
 and the notion of Julia set can be extended to any mapping on a Riemann surface to itself.
 For that we first need to define the Fatou set.
 
 The \emph{Fatou set} of a holomorphic map $ f $ is the set of points $ z $
 which have a neighborhood on which the 
 family $ \sequence{ f^n }_{ n \geq 0 } $
 of iterates of $ f $ forms a normal family.
In other words
 the point $ z $ belongs to the Fatou set of $ f $ if and only if
from  any subsequence of $ \sequence{ f^n }_{ n \geq 0 } $
 one can extract a (sub-)subsequence
 converging on
 a neighborhood of $ z $
 for the topology of local uniform convergence.
 The Fatou set is a totally invariant open set and the \emph{Julia set}
 is defined as the complement of the Fatou set in the domain of $ f $.
This means that the dynamics on the Fatou is stable while the Julia set contains the chaotic part of the dynamics.

When the Riemann surface in question is $ \complexnumbers $ the set of holomorphic functions is the set of entire functions, including the polynomials.
 On the Riemann sphere $ \riemannsphere $, the holomorphic functions are the \emph{rational maps}.
 Polynomial maps are also rational maps.
 A polynomial map is
 a rational map having a fixed point (identified with $ \infty \in \riemannsphere $) with no other preimage than itself.
 In particular $ \infty $ is a superattracting fixed point (see below) for any polynomial.

We will also need the following definitions.
A \emph{periodic point} for $ f $
 is a point $ z $ such that there exists $ p \geq 1 $ with $ f^p (z) = z $.
 The minimal value of $ p $ such that the above is satisfied is called the \emph{period} of $ z $.
 When the period is $ p = 1 $ a periodic point is simply called a \emph{fixed point}.
 When a point has a finite orbit but is not periodic it is called \emph{preperiodic}.
 
 The derivative of $ f^p $ at a periodic point of period $ p $ is called the \emph{multiplier}
 of the periodic point.
 The multiplier determines the local dynamics of $ f^p $ near the periodic point.
 Let $ \lambda $ be the multiplier of a periodic point $ z $. We have the following classification:
 \begin{enumerate}
  \item If $ \lambda = 0 $, the periodic point is called \emph{superattracting}.
  \item If $ \abs{ \lambda } < 1 $, the periodic point is called \emph{attracting} (superattracting is a special case of attracting).
  \item If $ \abs{ \lambda } = 1 $, the periodic point is called \emph{neutral}.
  \item If $ \abs{ \lambda } > 1 $, the periodic point is called \emph{repelling}.
 \end{enumerate}
In the first two cases the point $ z $ belongs to the Fatou set and has a \emph{basin of attraction}.
The basin of attraction 
 is an open neighborhood of $ z $ consisting of all of the points whose orbit under $ f^p $ converges to a point in the (finite) orbit of $ z $.

 The repelling periodic points belong to the Julia set and the Julia set is equal to the closure of the set of repelling periodic points of $ f $.
The neutral case is the most complicated (and interesting) and the point $ z $ might or might not belong to $ \juliaset (f) $
depending
on the map and, more importantly,
on the arithmetic properties of $ \lambda $.

\section{Complex dynamics and Kleinian groups}

The earliest picture of the Mandelbrot set\footnote{More precisely the conjectured interior of the Mandelbrot set.} to appear came from the study of discrete subgroups of M{\"o}bius transformations.
J{\o}rgensen \cite{Jorgensen1977} showed that a non elementary
\footnote{A subgroup of $\speciallinear (2, \complexnumbers) $ is \emph{elementary} if any pair of elements of infinite order 
have a common fixed point. 
Discreteness of elementary groups
can be checked in a easier way than non elementary groups.
}
subgroup of $ \speciallinear (2, \complexnumbers) $
is discrete if and only if all of its subgroups that are generated by two elements are discrete.
This result follows from an inequality
that J{\o}rgensen proved in an earlier work
\cite{Jorgensen1976}. 
This inequality known as \emph{J{\o}rgensen's inequality} is
a necessary condition for a group with two generators
to be discrete in $ \speciallinear (2, \complexnumbers) $.
 The proof consists of a rather simple argument by contradiction.
 Assuming that the group is not discrete one can easily find a pair of elements for which J{\o}rgensen's inequality is not satisfied.

The above results motivated the search for properties
 of subgroups of $ \projectivespeciallinear (2, \complexnumbers) $ generated by two elements that would imply discreteness.
 Brooks and Matelski \cite{BrooksMatelski1978}
 generalized J{\o}rgensen's result.
 This result can be stated as follows.
 Recall that an element $ \gamma $ of $ \projectivespeciallinear (2, \complexnumbers) $ is called \emph{loxodromic}
 if it is conjugated to $ z \mapsto k z $ for some $ k \in \complexnumbers \setcomplement \setof{ 0 } $ with $ \abs{ k } \neq 1 $
 or equivalently, if its squared trace $ \tracesquared \gamma $ is a complex number outside the closed interval $ \closedinterval{ 0, 4 } $.
\begin{thm}[Brooks and Matelski, 1978]
 Let $ \gamma_0, \gamma_1 $ be elements of $ \projectivespeciallinear (2, \complexnumbers) $ 
 with $ \gamma_0 $ loxodromic.
 Then there exists $ c = c (\gamma_0, \gamma_1) \in \complexnumbers $ and 
  $ z_0 = z_0 (\gamma_0, \gamma_1) \in \complexnumbers $
  such that if the subgroup generated by $ \gamma_0 $ and $ \gamma_1 $ is Kleinian
  then the set $ \setof{ z_n : n \in \positiveintegers } $ is discrete in $ \complexnumbers $,
  where the sequence $ z_n $ is defined by the following induction:
\begin{equation}
 z_{ n + 1 } = z_n^2 + c.
\end{equation}

The constants $ c $ and $ z_0 $ can be computed explicitly from $ \gamma_0 $ and $ \gamma_1 $.
\end{thm}

More precisely,
let
 $ \gamma_0, \gamma_1 \in \projectivespeciallinear (2, \complexnumbers) $
 with 
 $ \gamma_0 $ loxodromic.
Let $ \tau $ be the \emph{complex translation length} of  $ \gamma_0 $.
It is defined by the identity $ \tracesquared \gamma_0 = 4 \p{ \cosh (\tau / 2) }^2 $
with the normalizations $ \realpart \tau \geq 0 $ and $ \imaginarypart \tau \in \openclosedinterval{ -\pi, \pi } $
\footnote{
Equivalently $ \gamma_0 $ is conjugated to the map $ z \mapsto k z $ with $ k = e^{ \tau } $ and $ \realpart \tau \geq 0 $.}.
Then
\begin{equation}
 c = (1 - \cosh (\tau)) \cosh (\tau).
\end{equation}
Now define
for $ n \in  \integers_{ \geq 1 } $,
 $ \gamma_{ n + 1 } = \gamma_i \gamma_0 \gamma_n^{ -1 } $.
 Note that for $ n \geq 2 $, $ \gamma_n $ is loxodromic.
Let $ \delta_n $ be the \emph{complex distance} between the axis of $ \gamma_0 $ and the axis of $ \gamma_n $.
This complex number satisfies
 $ \realpart \delta_n \geq 0 $ and
 if we denote the fixed points of 
 $ \gamma_n $ by $ \alpha_n, \beta_n $
 (in order such that $ \alpha_n $ is repelling and $ \beta_n $ is attracting)
 then
 $ \p{ \cosh (\delta_n / 2) }^2 $
 is equal to the cross ratio of $ \alpha_0, \alpha_n, \beta_0, \beta_n $.
Then
\begin{equation}
  z_n = (1 - \cosh (\tau)) \cosh \delta_n.
\footnote{
An equivalent formulation is as follows.
Let $ \sigma $ be the squared trace of $ \gamma_0 $
and $ R $ the value of the cross ratio of
$ \alpha_0, \alpha_1, \beta_0, \beta_1 $. 
Then
\begin{equation}
 c = \p{ 2 - \sigma / 2 } \p{ \sigma / 2 - 1 } 
\end{equation}
and
\begin{equation}
 z_1 = \p{ 2 - \sigma / 2 } (2 R - 1). 
\end{equation}
Then $ \frac{ z_n }{ 2 - \sigma / 2 } $
is equal to the image of
the complex distance between the axis of $ \gamma_n $ and the axis of $ \gamma_0 $
by the function $ \cosh $.
}
\end{equation}
 
 If the group generated by $ \gamma_0 $ and $ \gamma_1 $ is discrete then the set
$ \setof{ \cosh (\delta_n) : n \in \integers_{ \geq 2 } } $
is discrete in $ \complexnumbers $.
The theorem follows from the following inductive relation on the sequence of $ \delta_n $  \cite[p.67]{BrooksMatelski1978}:
\begin{equation}
 \cosh (\delta_{ n + 1 }) = (1 - \cosh (\tau)) \p{ \cosh (\delta_n) }^2 + \cosh \tau.
\end{equation}
In their article they proceed to draw the filled Julia set of $ z^2 + 0.1 + 0.6 i $ and
the set of $ c \in \complexnumbers $ for which $ z \mapsto z^2 + c $
has a stable periodic orbit.
It is noteworthy that one of the important features of the field of complex dynamics
 in the years following the work of Brooks and Matelski is the use of computer graphics in an exploratory way.

Works on the question of the discreteness of groups of M{\"o}bius transformations generated by certain generators,
 and in particular on generalizations of J{\/o}rgensen's inequalities, neither started nor ended with the above example
 (e.g. \cite{Jorgensen1976,Keen1983,Tan1989,KeenSeries1994,Cao1995,Bowditch1998,Klimenko2001,KlimenkoKopteva2002}).
 One of the
 important developments appears in the work of Gehring and Martin \cite{GehringMartin1989}.
 Their main theorem is as follows.
 
 \begin{thm}[Gehring and Martin, 1989]
  Assume that the group generated by $ \gamma_0 \in \projectivespeciallinear (2, \complexnumbers) $ and $ \gamma_1 \in \projectivespeciallinear (2, \complexnumbers) $
  is Kleinian and $ \gamma_0 $ is loxodromic.
  Define 
  \begin{equation}
   z_0 = \operatorname{tr}  \commutator{ \gamma_0, \gamma_1 }  - 2
  \end{equation}
  and
  \begin{equation}
   \beta = \tracesquared \gamma_0 - 4,
  \end{equation}
  where $ \commutator{ \gamma_0, \gamma_1 } = \gamma_0 \gamma_1 \gamma_0^{ -1 } \gamma_1^{ -1 } $ is the commutator.
Let $ K (P_{ \beta }) $ be the filled Julia set of $ z \mapsto P_{ \beta } (z) = z^2 - \beta z $.

 Then either $ z_0 \notin K (P_{ \beta }) $
 or $ z_0 $ is preperiodic under the iteration of $ P_{ \beta } $.
 If $ z_0 $ is preperiodic, its orbit never lands on the fixed point $ 0 $.
 
 Moreover if $ z_0 $ is preperiodic there are nontrivial conjugacy relations between $ \gamma_0 $ and  $ \gamma_1 $.
 \end{thm}

These results can be related to other types of link that have been established
between the iteration of rational maps and Kleinian groups.
 These include famously Sullivan's dictionary (compare \cite{Sullivan1983,Sullivan1985,McMullen1991,McMullen1994})
 but
 also the study of the coupling of the dynamics of rational maps with M{\"obius} transformations through a procedure called \emph{mating} (see e.g. 
 \cite{BullettPenrose1994,BullettPenrose1999,BullettHaissinsky2007, BullettLomonaco2016ax}).

\section{Newton's method and other numerical methods}

Since for most polynomials there is no simple formula that expresses the roots in terms of the coefficients
 one has to use iterative methods to find numerical approximations of their roots.
 A classic method is Newton's iterative scheme.
 It is based on the idea that the function whose roots are to be found can be locally replaced by its first order approximation.
 An approximation of a root is inductively computed using this local approximation of the function.
 In precise terms, if we want to solve the equation
 \begin{equation}
  P (z) = 0
 \end{equation}
we define a sequence of approximations of some root by picking a guess $ z_0 $ and then defining the sequence
\begin{equation}
 z_{ n + 1 } = N (z_n)
\end{equation}
where $ N_P $ is the \emph{Newton map} and is defined as
\begin{equation}
 N_P (z) = z - \frac{ P (z) }{ P' (z) }.
\end{equation}
For a large set of choices of $ z_0 $ the sequence $ \sequence{ z_n }_n $ will indeed converge to a root of $ P $.
Since the invention of Newton's method many other methods for finding the roots of polynomials have been found.
Despite its simplicity Newton's method is already quite efficient.
This simplicity and its old age has allowed a good understanding of the dynamics of Newton's method.
When applied to complex analytic equations, in particular when $ P $ is a polynomial, Newton's method can be studied by using complex dynamics.

The study of the dynamics of the Newton map is an old (for example, \cite{Cayley1879}) and rich topic. 
Here we are only exploring a very small portion of the theory.
In particular we only look at methods for finding roots of a polynomial.
The study of Newton's method in complex dynamics is not restricted to this case,
 see for example \cite{Haruta1992,Haruta1999,Berweiler1993} and more recently \cite{BaranskiFagellaJarqueKarpinska2018}.

We will focus on two aspects of the theoretical study of Newton's method and related root finding algorithms.
 Firstly we will ask the question of
 how big is the set of pairs map-and-initial-guess $ (f, z_0) $ for which the method converges.
 Then we will look for an algorithm to find all the roots of a given polynomial.

\subsection{Genericity of convergence}

Given an iterative algorithm it is natural to ask for which initial values this algorithm will converge.
One can also ask for which function we are guaranteed to find the roots by using the algorithm.
In the best case the algorithm would converge to all or almost all (in the sense of measure)
 pairs $ (P, z) \in \polynomials (d) \times \complexnumbers $ of polynomials and initial guesses.
In the context of Newton's method and more general root finding algorithms
 this idea as been conceptualized by Smale.
 Smale introduced in \cite{Smale1987}
 the notion of \emph{generally convergent purely iterative algorithm} (GCPIA).
 
 A \emph{purely iterative algorithm} is given by a map $ T (P, z) = T_P (z) $
 which depends rationally on $ z $ and on the coefficients of $ P $.
 A purely iterative algorithm is \emph{generally convergent} if
   there exists a dense open set $ \Omega \subset \polynomials (d) \times \complexnumbers $
   of full measure
   such that
 for
 all
 $ (P, z) \in \Omega $
  the sequence $ \sequence{ T_P^n (z) }_n $
  converges to a root of $ P $.
  A \emph{GPCIA} is a purely iterative algorithm which is generally convergent.
  
 An alternative definition requires only that $ \Omega $ is open and dense but not necessarily of full measure \cite{McMullen1988}.
  To distinguish them from GCPIA we will call such algorithms
 \emph{GCPIAM}.

Newton's method is not a GCPIA for $ \polynomials (d) $, $ d > 2 $. Indeed there are many examples of polynomials for which the Newton map has other attracting basins than
 the ones of the roots, see for example \cite{Hurley1986}.
Using deep results in complex dynamics McMullen was able to
 show that there is no GCPIA for $ \polynomials (d) $ with $ d \geq 4 $
 and gave a complete classification of GCPIA for $ d = 2, 3 $.

Before stating McMullen's result we need the following definition.
The \emph{centralizer} $ C (T) $ of a rational map $ T $ is defined as
 the subgroup of M{\"o}bius transformations which commute with $ T $.

\begin{thm}[McMullen, \cite{McMullen1987}, Theorem 1.1] \label{thm: GCPIA for polynomials}
\hspace{2em} 
 \begin{enumerate}
  \item There is no GCPIAM for $ \polynomials (d) $ with $ d \geq 4 $.
  
  \item Let $ T $ be a purely iterative algorithm defined over $ \polynomials (3) $.
  Then $ T $ is a GCPIAM if and only if there exists a rational map $ T_0 : \riemannsphere \to \riemannsphere $ such that the following are true.
  \begin{enumerate}
   \item There is $ U_0 \subset \complexnumbers $ open and dense in $ \complexnumbers $ such that
    for all $ z \in U_0 $, $ T_0^n (z) $ converges to a root of $ P_0 (z) = z^3 - 1 $. 
   \item The centralizer $ C (T) $ contains the group of M{\"o}bius transformations
   permuting the roots of $ P_0 $.
   \item For all $ P \in \polynomials (3) $ with no multiple root, $ T_P = M_P \composedwith T_0 \composedwith M_P^{ -1 } $
   where $ M_P $ is a M{\"o}bius transformation mapping the roots of $ P_0 $ to the roots of $ P $.
  \end{enumerate}
   
   \item Let $ T $ be a purely iterative algorithm defined over $ \polynomials (2) $.
  Then $ T $ is a GCPIAM if and only if there exists a rational map $ T_0 : \riemannsphere \to \riemannsphere $ 
  and a rational function $ M : \polynomials (2) \to \projectivespeciallinear (2, \complexnumbers) / C (T_0) $
  such that the following are true.
  \begin{enumerate}
   \item There is $ U_0 \subset \complexnumbers $ open and dense in $ \complexnumbers $ such that
    for all $ z \in U_0 $, $ T_0^n (z) $ converges to a root of $ P_0 (z) = z^2 - 1 $.
   \item The centralizer of $ T_0 $ contains $ z \mapsto - z $.
   \item If $ P $ has no multiple roots then $ M (P) $ maps the roots of $  P_0 $ to the roots of $ P $
   and $ T_P = M_P \composedwith T_0 \composedwith M_P^{ -1 } $ for some representant $ M_P $ of $ M (P) $ in $ \projectivespeciallinear (2, \complexnumbers) $.
  \end{enumerate}
  
 \end{enumerate}

\end{thm}

In particular the following examples are GCPIA (see \cite{McMullen1987}, Proposition 1.2):
\begin{enumerate}
 \item Newton's method for quadratic polynomials.
 \item The Newton's map of the rational map 
 \begin{equation}\label{eq: cubic universal GCPIA map}
  f (z) = \frac{ z^3 + a z + b }{ 3 a z^2 + 9 b z - a^2 }
 \end{equation}
  is a GCPIA for the cubic polynomials of the form
  \begin{equation} \label{eq: cubic reduced form}
   P (z) = z^3 + a z + b.
  \end{equation}
\end{enumerate}
The above are also characterized by their fast convergence due to the fact that the roots are superattracting fixed points of the map $ T_P $
(compare above reference).
Note that finding the roots of a cubic polynomial can easily be replaced by the problem of finding the roots of some $ P $
in the form \eqref{eq: cubic reduced form}. 
The Newton map for \eqref{eq: cubic universal GCPIA map} is \emph{expanding} and its only Fatou components are in the basin of some root of $ P $.
This implies that it is not just a GCPIAM but also a GCPIA as stated.

When the condition on the complex analycity of the mapping $ T $
 is relaxed into real analycity (that is by allowing complex conjugate in the formulas), a GCPIA exists for any degree
 \cite{ShubSmale1986}.

McMullen later refined their results in \cite{McMullen1988}. 
 In that article they explain that \emph{braiding} of the roots when going around a polynomial with multiple roots
 prevents the existence of a mapping $ T_P $
 which is a GCPIA for the polynomials of degree $ d \geq 4 $.
 This is also a very interesting article where complex dynamics is used for studying many aspects of the GCPIAs.

The proof of Theorem \ref{thm: GCPIA for polynomials} relies on deep results of complex dynamics.
These include the celebrated work of Ma{\~n}{\'e}, Sad and Sullivan \cite{ManeSadSullivan1983}
and Thurston's work on the characterization of \emph{postcritically finite rational maps}.
A postcritically finite rational map is a rational map $ f : \riemannsphere \to \riemannsphere $
 such that all of its critical points are either periodic or preperiodic.
 
If $ T $ is a GCPIAM for polynomials of degree $ d \geq 2 $
 then $ \family{ T_P }{ P \in \polynomials (d) } $
  forms a \emph{stable algebraic family}.
This means that this  
  is a family of rational maps of fixed degree
  depending rationally on the coefficients of $ P $
  (\emph{algebraic family})
  and there is a uniform bound on the periods of attracting cycle (\emph{stable}).
Indeed the roots are the only attracting periodic points of the family for generic points
 and there is no bifurcation in the sense of \cite{ManeSadSullivan1983}.

From a result of Thurston it follows that 
 stable algebraic families
 either are \emph{trivial}
 (all the elements in the family are conjugated to each other by a M{\"o}bius transformation)
 or consist of Latt{\`e}s examples (see \cite{Lattes1918} or \cite{Milnor2006book}, Definition 7.4 for a definition).
 The latter case is excluded for a GCPIAM since
 the Julia set of a Latt{\`e}s example consists of the whole Riemann sphere.
It follows from the rigidity of M{\"o}bius transformations that a GCPIA cannot exists for $ d \geq 4 $.

To get around the problem of the non existence of GCPIA
 one can consider instead 
 \emph{towers of algorithms} as defined in \cite{DoyleMcMullen1989}.
Then it can be shown that the
roots of a polynomial of degree $ d $ can be computed by a general tower of algorithms
 if and only if $ d \leq 5 $ (\cite{DoyleMcMullen1989}, Corollary 4.3).
An explicit algorithm for the quintic in given in the appendix of \cite{DoyleMcMullen1989}.

Crass
\cite{Crass2001}, \cite{Crass2002}
has provided methods for solving the quintics and equations of higher degree 
 in a similar manner.
 These methods involve the iteration of holomorphic maps in higher dimensional complex projective spaces.
Subsequent developments also include
\cite{Crass2014}.

\subsection{Finding all the roots}

One of the remarkable feature of the theory is that it can be used to describe an explicit strategy for finding all the roots of a polynomial with certainty.

Early works on the maximal complexity of Newton's method applied to the search of roots of complex polynomials include Manning's \cite{Manning1992}.
 This works contains the description of an implementable algorithm that ensures the finding of at least one root with complexity bounded a priori by a constant depending
 only on the degree $ d $ (note that their result applies only for $ d \geq 10 $).
 This is based on the fact that the Newton map $ N_P $ has a repelling fixed point on the Riemann sphere at $ \infty $,
 explicit bounds on the behavior of $ N_P $
 and distortions estimates coming from complex analysis.

The question of finding a choice of initial guesses that would guarantee
 finding all the roots of the polynomial was answered by Hubbard, Schleicher and Sutherland in \cite{HubbardSchleicherSutherland2001}.
 The set they produce depends only on the degree $ d $ of $ P $.

 Let $ \mathcal{P}_d $ be the set of polynomials of degree $ d $ with all the roots inside the open unit disk $ \unitdisk $.
Note that there is a simple method to substitute the problem of finding all the roots of an element of $ \mathcal{ P }_d $
for the problem of finding roots of some arbitrary polynomial. 

\begin{thm}[\cite{HubbardSchleicherSutherland2001}] \label{thm: finding all the roots using NM}
 Let $ d \geq 2 $.
 There exists $ S_d \subset \complexnumbers $ finite with at most $ 1.11 d \p{ \log d }^2 $ elements
 such that for all $ P \in \mathcal{P}_d $ and all root $ \xi $ of $ P $
 there exists $ s \in S_d $
 such that
 $ N_P^n (s) \tendsto \xi $ as $ n \tendsto \infty $.
\end{thm}

Let $ \eps > 0 $. By compactness there exists $ n = n (d, \eps) $
such that for all $ P \in \mathcal{P}_d $  and all $ \xi $ root of $ P $
there exists $ s \in S_d $ such that
\begin{equation}
 \abs{ N_P^n (s) - \xi } \leq \eps.
\end{equation}
This ensures that the algorithm effectively finds all the roots in finite time.
Schleicher's article \cite{Schleicher2002}
provides explicit estimates on $ n (d, \eps) $.
In theory each guess could require a large number of iterations as the degree becomes large.

The article \cite{HubbardSchleicherSutherland2001} also contains an explicit construction for 
the set $ S_d $ and finer and better results for when the polynomial is real.
The authors use their own algorithm to compute approximations to the invariant measure of H{\'e}non mappings.

In \cite{SchleicherStoll2017}
Schleicher and Stoll
give
a slightly different version of the algorithm mentioned above.
 This reference also contains many remarks on the implementation and possible improvements.
 Using some numerical experiments they checked
 that the theoretically possible large number of iterations
 (larger than $ d^2 $ with $ d \approx 10^6 $)
 was not a problem in practice for the specific problems they were looking at.
 They used it to find the centers of hyperbolic components of the Mandelbrot set
 and periodic points of iterated polynomials.

 The roots are attracting fixed points of the Newton map $ N_P $.
The \emph{basin of a root} is the set of points whose orbit converge to the root under the iteration of
 $ N_P $. This is an open set.
 The \emph{immediate basin} of a root is the connected component of the basin containing the root.
The proof of Theorem \ref{thm: finding all the roots using NM}
 builds on previous results relating to the shape of the immediate basins of the roots 
 such as \cite{Manning1992} and \cite{Przytycki1989}.

A summary of the proof is as follows.
 The only fixed points of the Newton map $ N_P $
 are the roots of $ P $ and $ \infty $.
 The roots are either superattracting or attracting with multiplier $ 1 - 1/k $ for some integer $ k \geq 2 $.
 The fixed point at $ \infty $ is repelling.
 
 Let $ \sequence{ \xi_i }_i $ be the roots of $ P $ and let $ U_{ \xi_i } $ be the immediate basin of $ \xi_i $.
 Define also $ m_{ \xi_i } $ as the number of critical points of $ N_P $ (counted with multiplicity) inside $ U_{ \xi_i } $.
 From \cite{HubbardSchleicherSutherland2001}, Proposition 6, it follows that
 $ U_{ \xi_i } $ has $ m_{ \xi_i } $ accesses to $ \infty $
 \footnote{ That is the complement of some large disk in $ U_{ \xi_i } $ has $ m_{ \xi_i } $ components accumulating to $ \infty $. }.
The idea is to constrain  the geometry of these accesses.
 
 Pick a root $ \xi = \xi_i $.
From \cite{Przytycki1989} (see also \cite{Shishikura2009}) we know that $ U_{ \xi } $ is simply connected.
Let $ \ph : \unitdisk \to U_{ \xi } $ be a conformal isomorphism normalized so that $ \ph (0) = \xi $
and define
\begin{equation}
 f \defeq \ph^{ -1 } \composedwith N_P \composedwith \ph.
\end{equation}
The mapping $ f $ is proper of degree $ m + 1 $ where $ m = m_{ \xi } $.
This mapping can be extended by reflection into a rational map $ f : \riemannsphere \to \riemannsphere $
 of degree $ m + 1 $.
 
The rational map $ f $  has $ m + 2 $ fixed points (counted with multiplicity).
 The point $ 0 $ is a (super)attracting fixed point of $ f $.
 By symmetry this is also the case for the point $ \infty $.
 The respective multipliers $ \lambda_0, \lambda_{ m + 1 } $ of $ 0 $ and $ \infty $
 are either both equal to  $ 0 $ or to $ 1 - 1 / k $ for some integer $ k \geq 2 $.
 It also follows from the symmetry that the other fixed points
 $ \zeta_1, \dots, \zeta_m $
 lie on the unit circle
 and their respective multipliers are positive real numbers $ \lambda_j > 1 $.

The holomorphic fixed point formula applied to $ f $ (see e.g. \cite[Section 12]{Milnor2006book})
states that
\begin{equation}
 \familysum{ j = 0 }{ m + 1 } \frac{ 1 }{ \lambda_j - 1 } = - 1.
\end{equation}
Hence
\begin{equation}
 \familysum{ j = 1 }{ m } \frac{ 1 }{ \lambda_j - 1 } \geq 1. 
\end{equation}
It follows that there must be at least one $ j $ such that $ \lambda_j - 1 \leq m $.

The quotient of the corresponding channel for $ N_P $ by the dynamics of $ N_P $
 is an annulus of modulus $ \frac{ \pi }{ \log \lambda_j } $.
 Indeed this is the value of the modulus of the annulus obtained by
 taking the quotient of the upper plane by the action of $ z \mapsto \lambda_j z $.
 Since the degree of $ f $ is at most equal to the degree of $ N_P $
 it follows that
 \begin{equation}
  \frac{ \pi }{ \log \lambda_j } \geq \frac{ \pi }{ \log (m + 1) } \geq \frac{ \pi }{ \log d }.
 \end{equation}
Having such a lower bound on the modulus allows to find places where the channel must have a definite extent.
 This is made precise in  \cite[Section 5]{HubbardSchleicherSutherland2001}.
Using this, one can pick points independently of $ P $ such that at least one of them is in $ U_{ \xi } $.

%
%
%
%
%
%
%
%
%
%
%
%
%
%

%
%

%

\section{Hierarchical Ising and Potts models}

The Ising and Potts models are mathematical models from solid state physics.
 The Ising model relates to the ferromagnetic properties of a material.
 At the base of both models lies a graph whose vertexes represent the locus of a particle/atom and the edges the interaction between these particles.
 Each vertex is characterized by a state chosen among a finite set of possible values.
 For the Ising model this set has $ 2 $ element while for the Potts model
 the number of possible states is some positive integer $ q \geq 2 $.
 
 The Hamiltonian of the system can be computed explicitly for any state.
  The temperature $ T $,
  interaction constant $ J $ and the (possibly $ 0 $) magnetic field $ h $
  appear as parameters in the Hamiltonian.
  From the Hamiltonian one can derive a formula for the partition function.

\begin{figure}
 \centering
  \def\svgwidth{\columnwidth}
  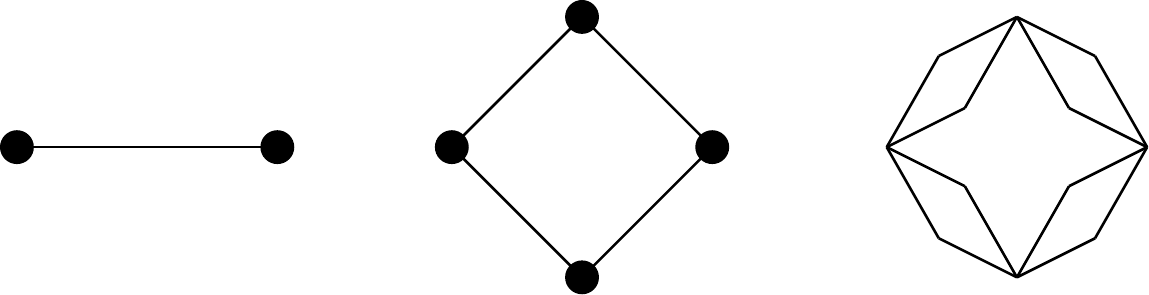
  \caption{The first levels of the diamond lattice hierarchy ($ b = 2 $).} \label{fig: hierarchical model lattices}
 \end{figure}

\begin{figure} 
 \centering
  \def\svgwidth{\columnwidth}
  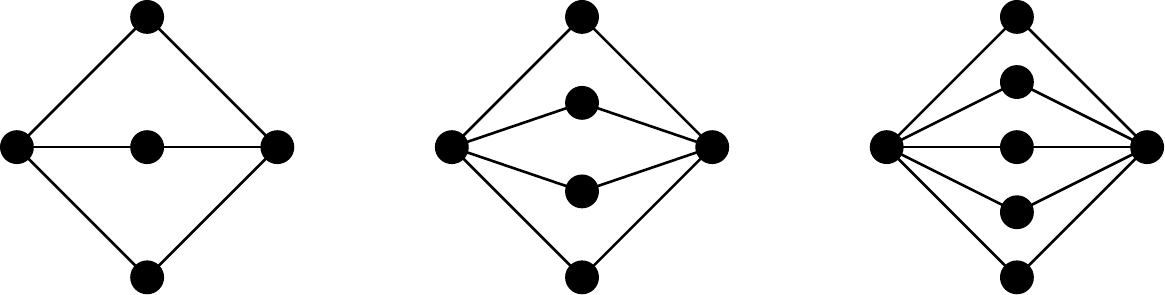
  \caption{Examples of diamond lattices with $ b = 3, 4, 5 $.} \label{fig: hierarchical model lattices b gt 2}
 \end{figure}

 A \emph{hierarchical lattice} consists of a refining sequence of finite graphs on which the Hamiltonian is computed successively.
  For example a \emph{diamond hierarchical lattice} can be defined as follows.
  The first graph consists of a pair of vertexes joined by a single edge.
  The refining consists in replacing each edge by two pairs of edges each connecting one of the previous vertexes to a new vertex in the middle.
  The refining procedure is illustrated in Figure \ref{fig: hierarchical model lattices}.
  Instead of replacing each edges by $ 2 $ branches (pairs of edges), one could also insert $ b \geq 2 $ branches,
   see Figure \ref{fig: hierarchical model lattices b gt 2}.
   The integer $ b $ is the parameter characterizing a diamond hierarchical lattice.

 The passage from the partition function of one level of the hierarchy to the next level is performed by a \emph{renormalization group transformation}.
 This transformation depends on the variable $ J $.
 For a diamond hierarchical lattice model the renormalization group transformation is identified as a rational map $ f : z \mapsto f (z) $, where $ z = z (J, T) $,
 and its dynamics has a physical relevance.
 For example the zeros of the partition functions in the thermodynamic limit
 (i.e. when the level in the hierarchy tends to $ \infty $)
 converges to the Julia set of $ f $.
 The dynamics of this map is the focus of the study of hierarchical Potts/Ising models in complex dynamics.
 
These models generally do not really represent actual physical systems
  but are instead used to try understand what type of properties
  more complicated and realistic model could have.
  In general one cannot hope to have an explicit formula for the renormalization transformation of a realistic model.
 
 Hierarchical models are described in \cite{McKayBerkerKirkpatrick1982}, \cite{BerkerMcKay1984} and \cite{DerridaDeSezeItzykson1983}.
 In the latter, Derrida, De Seze and Itzykson study the $ q $-state Potts model on a diamond hierarchical lattice with $ b = 2 $.
The renormalization map $ f $ for this model can be computed explicitely:
 \begin{equation}
  f (z) = \p{ \frac{ z^2 + q - 1 }{ 2 z + q - 2 } }^2.
 \end{equation}
They provide several pictures of the Julia sets corresponding to different values of $ q $ in an attempt to get an idea of their fractal structure.
This work has been followed by \cite{DerridaItzyksonLuck1984} where the geometric properties of the Julia set of $ f $
 are used to extract informations about the model. Those are mainly numerical studies.

Another type of hierarchical model is presented in \cite{BoscoRosa1987}.
 The renormalization transformation can also be identified to a rational map.
 The authors study the structure and Hausdorff dimension of the corresponding Julia set.
 For another model see also \cite{AnanikianGhulghazaryan2000}, \cite{GhulghazaryanAnanikyanJonassen2003}.

%

An important occurence of the utilization of complex dynamics to study the Ising model is the work of
Bleher and Lyubich
\cite{BleherLyubich1991}.
They study the Ising model on the diamond hierarchical lattice for arbitrary values of $ b \geq 2 $ (see also \cite{BleherZalis1989}).
In that case the renormalization transformation is represented by the rational map
\begin{equation}
 f (z) = \frac{ 4 z^b }{ \p{ 1 + z^b }^2 }.
\end{equation}
The points $ 0 $ and $ 1 $ are superattracting fixed points of $ f $.
Denote the immediate basin of $ 0 $ by $ \Omega_0 $.
The free energy can be expressed as
\begin{equation}
 F (z) = \familysum{ n = 0 }{ \infty } \frac{ 1 }{ \p{ 2 b }^n } g \composedwith f^n (z)
\end{equation}
where $ g (z) = \log ( 1 + z^b ) $.
Bleher and Lyubich showed that $ F $ is analytic on $ \Omega_0 $
and that the boundary of $ \Omega_0 $ is a natural boundary of analyticity for this function
(i.e. analytic continuation is not possible along any path that crosses $ \boundary \Omega_0 $).
They also derive some physically relevant properties of the model.

The Fatou set of $ f $ consists of the respective basins of attractions of $ 0 $ and $ 1 $.
Before proceeding to study the properties of the free energy $ F $,
 they first showed that $ \Omega_0 $ is a Jordan domain.
 Recall that a  \emph{quasicircle} is the image of a round circle by a quasiconformal homeomorphism of $ \complexnumbers $.
 Bleher and Lyubich showed that the boundary of
 $ \Omega_0 $ is a quasicircle.
Although this proof is rather simple
 it exemplifies the use of a powerful tool of complex dynamics:
 the theory of polynomial like maps \cite{DouadyHubbard1985polylike}.
 This theory explains why copies of the Mandelbrot set seem to appear in every parameter space of holomorphic dynamical systems.

Let $ d \geq 2 $.
A \emph{polynomial like map of degree $ d $} is a triplet $ \sequence{ U, U', f } $
 where $ U $ and $ U' $ are simply connected open subsets of $ \complexnumbers $
 such that $ U' $ is compactly contained in $ U $
 and $ f : U' \to U $ is a proper holomorphic mapping of degree $ d $.
The \emph{filled Julia set} $ \filledjuliaset (f) $ of a polynomial like mapping is the set of points whose orbit stays inside the domain $ U' $,
\begin{equation}
 K (f) \defeq \setof{ z \in U' : \forall n, f^n (z) \in U'}.
\end{equation}

The relevance of polynomial like mappings derives from Douady and Hubbard's straightening theorem.
\begin{thm}[Douady, Hubbard, \cite{DouadyHubbard1985polylike}, Theorem 1]
 Let $ (U, U', f) $ be a polynomial like mapping of degree $ d \geq 2 $.
 Then there exists a quasiconformal map $ \psi : \complexnumbers \to \complexnumbers $
 and a polynomial $ P $ of degree $ d $
 such that 
 \begin{equation}
  \psi \composedwith f = g \composedwith \psi
 \end{equation}
 on some neighborhood of $ K (f) $
 and $ \dbar \psi = 0 $ almost everywhere on $ K (f) $.

 Moreover if $ K (f) $ is connected, then the polynomial $ P $ is unique up to conjugation by an affine map.
\end{thm}

Thanks to a fine analysis of the map $ f $,
Bleher and Lyubich showed that $ f $ is polynomial like of degree $ b $ on a neighborhood of the closure of $ \Omega_0 $.
Since it has a superattracting fixed point of degree $ b $ at $ 0 $,
the straightening of $ f $ is conjugated to the polynomial $ z \mapsto z^b $.
Since $ \Omega_0 $ is the basin of $ 0 $,
it follows that $ \boundary \Omega_0 $ is the image of the circle $ \setof{ z : \abs{ z } = 1 } $
by a quasiconformal map of the plane, hence it is a quasicircle.

The use of complex dynamics in the field has continued after this work, for example in \cite{BleherLyubichRoeder2017}.

\section{Other connections}  \label{sec: other}

There are many other applications of complex dynamics.
 Eremenko has mentioned other connections in a talk \cite{Eremenko2018talk}
  about the interaction between function theory and complex dynamics.
  
A surprising application is related to gravitational lensing.
 In \cite{KhavinsonSwiatek2003} Kahvinson and {\'S}wi{\k{a}}tek
 solved the Sheil-Small and Wilmhurst conjecture.
 This states that if $ P $ is a polynomial of degree $ n \geq 2 $
 then the harmonic polynomial $ z - \complexconjugate{ P (z) } $
 has at most $ 3 n - 2 $ zeros.
 
 Their proof relies on the classical fact from complex dynamics that
  any attracting or parabolic periodic point attracts at least one critical point.
  This can be applied to the holomorphic polynomial $ Q (z) = \complexconjugate{ P \p{ \complexconjugate{ P (z) } } } $.

  It turns out that this solution has an application in astrophysics exposed in \cite{KhavinsonNeumann2006}.
  A similar argument can be used when one replaces the polynomial $ P $ by a rational function $ R $.
  This gives the following theorem.
\begin{thm}[Khavinson, Neumann, \cite{KhavinsonNeumann2006}]
 Let $ R $ be a rational function
of degree $ n \geq 2 $.
Then the equation
$ z = \complexconjugate{ R (z) } $
has at most $ 5 n - 5 $ solutions.
\end{thm}
In astrophysics the lensing effect produced by the gravity coming from $ n $ point like objects can be modelled
 via a lens equation (see for example \cite{Witt1990} and \cite{Rhie2001ax}).
A corollary (\cite{KhavinsonNeumann2006}, Corollary 1) of the above theorem gives an explicit upper bound on the number of images that such model can produce.
  See \cite{KhavinsonNeumann2008} for further developments.
  For more details about this the reader is advised to consult the excellent \cite{Roeder2016}.


\bibliographystyle{plain}
\bibliography{connections-of-dh}

\end{document}

%% file: math-lib.tex
\usepackage{bbm}

\usepackage{calrsfs}

\usepackage{color}

\usepackage{hyperref}

\usepackage{enumitem}
\setenumerate[1]{label={\arabic*.}}

\newcommand{\hide}[1]{}

\newcommand{\eps}{\varepsilon}
\renewcommand{\epsilon}{\eps}
\newcommand{\ph}{\varphi}
\renewcommand{\phi}{\ph}

\newcommand{\overbar}[1]{\overline{#1}}

\newcommand*{\defeq}{\mathrel{\vcenter{\baselineskip0.5ex \lineskiplimit0pt
                     \hbox{\scriptsize.}\hbox{\scriptsize.}}}%
                     =}

\newtheorem{thm}{Theorem}[section]

\newcommand{\setof}[1]{\left\{#1\right\}}

\newcommand{\setcomplement}{{\mathbin{\backslash}}}

\newcommand{\composedwith}{\mathbin{\circ}}

\newcommand{\family}[2]{\left(#1\right)_{#2}}

\newcommand{\complexnumbers}{{\mathbbm{C}}}
\newcommand{\integers}{{\mathbbm{Z}}}

\newcommand{\positiveintegers}{{\integers_{>0}}}

\newcommand{\absolutevalue}[1]{\left|#1\right|}

\newcommand{\prioritize}[1]{\left(#1\right)}

\newcommand{\sequence}[1]{\left(#1\right)}
\newcommand{\tendsto}{{\to}}
\newcommand{\topologicalboundary}{{\partial}}

\newcommand{\closedinterval}[1]{\left[#1\right]}

\newcommand{\openclosedinterval}[1]{\left]#1\right]}

\newcommand{\familysum}[2]{{\sum\limits_{#1}^{#2}}\ }

\newcommand{\complexconjugate}[1]{\overline{#1}}

\newcommand{\imaginarypart}{\operatorname{Im}}

\newcommand{\realpart}{\operatorname{Re}}
\newcommand{\riemannsphere}{\widehat{{\complexnumbers}}}

\newcommand{\unitdisk}{{\mathbb{D}}}

\newcommand{\filledjuliaset}{K}

\newcommand{\juliaset}{J}

\newcommand{\commutator}[1]{\left[#1\right]}

\newcommand{\abs}[1]{\absolutevalue{#1}}

\newcommand{\boundary}{\topologicalboundary}

\newcommand{\p}{\prioritize}



%% file: diamond-lattice.pdf_tex
\begingroup%
  \makeatletter%
  \providecommand\color[2][]{%
    \errmessage{(Inkscape) Color is used for the text in Inkscape, but the package 'color.sty' is not loaded}%
    \renewcommand\color[2][]{}%
  }%
  \providecommand\transparent[1]{%
    \errmessage{(Inkscape) Transparency is used (non-zero) for the text in Inkscape, but the package 'transparent.sty' is not loaded}%
    \renewcommand\transparent[1]{}%
  }%
  \providecommand\rotatebox[2]{#2}%
  \newcommand*\fsize{\dimexpr\f@size pt\relax}%
  \newcommand*\lineheight[1]{\fontsize{\fsize}{#1\fsize}\selectfont}%
  \ifx\svgwidth\undefined%
    \setlength{\unitlength}{330.70059118bp}%
    \ifx\svgscale\undefined%
      \relax%
    \else%
      \setlength{\unitlength}{\unitlength * \real{\svgscale}}%
    \fi%
  \else%
    \setlength{\unitlength}{\svgwidth}%
  \fi%
  \global\let\svgwidth\undefined%
  \global\let\svgscale\undefined%
  \makeatother%
  \begin{picture}(1,0.25627411)%
    \lineheight{1}%
    \setlength\tabcolsep{0pt}%
    \put(0,0){\includegraphics[width=\unitlength,page=1]{diamond-lattice.pdf}}%
  \end{picture}%
\endgroup%

%% file: diamond-lattice-gt-2.pdf_tex
\begingroup%
  \makeatletter%
  \providecommand\color[2][]{%
    \errmessage{(Inkscape) Color is used for the text in Inkscape, but the package 'color.sty' is not loaded}%
    \renewcommand\color[2][]{}%
  }%
  \providecommand\transparent[1]{%
    \errmessage{(Inkscape) Transparency is used (non-zero) for the text in Inkscape, but the package 'transparent.sty' is not loaded}%
    \renewcommand\transparent[1]{}%
  }%
  \providecommand\rotatebox[2]{#2}%
  \newcommand*\fsize{\dimexpr\f@size pt\relax}%
  \newcommand*\lineheight[1]{\fontsize{\fsize}{#1\fsize}\selectfont}%
  \ifx\svgwidth\undefined%
    \setlength{\unitlength}{335.25bp}%
    \ifx\svgscale\undefined%
      \relax%
    \else%
      \setlength{\unitlength}{\unitlength * \real{\svgscale}}%
    \fi%
  \else%
    \setlength{\unitlength}{\svgwidth}%
  \fi%
  \global\let\svgwidth\undefined%
  \global\let\svgscale\undefined%
  \makeatother%
  \begin{picture}(1,0.25279642)%
    \lineheight{1}%
    \setlength\tabcolsep{0pt}%
    \put(0,0){\includegraphics[width=\unitlength,page=1]{diamond-lattice-gt-2.pdf}}%
  \end{picture}%
\endgroup%

%% file: connections-of-dh.bbl
\begin{thebibliography}{10}

\bibitem{AnanikianGhulghazaryan2000}
N.~S. Ananikian and R.~G. Ghulghazaryan.
\newblock Yang-{L}ee and {F}isher zeros of multisite interaction {I}sing models
  on the {C}ayley-type lattices.
\newblock {\em Phys. Lett. A}, 277(4-5):249--256, 2000.

\bibitem{BaranskiFagellaJarqueKarpinska2018}
Krzysztof Bara\'{n}ski, N\'{u}ria Fagella, Xavier Jarque, and Bogus{\l{}}awa
  Karpi\'{n}ska.
\newblock Connectivity of {J}ulia sets of {N}ewton maps: a unified approach.
\newblock {\em Rev. Mat. Iberoam.}, 34(3):1211--1228, 2018.

\bibitem{Beardon1991}
Alan~F. Beardon.
\newblock {\em Iteration of rational functions}, volume 132 of {\em Graduate
  Texts in Mathematics}.
\newblock Springer-Verlag, New York, 1991.
\newblock Complex analytic dynamical systems.

\bibitem{Berweiler1993}
Walter Bergweiler.
\newblock Newton's method and a class of meromorphic functions without
  wandering domains.
\newblock {\em Ergodic Theory Dynam. Systems}, 13(2):231--247, 1993.

\bibitem{BerkerMcKay1984}
A.~Nihat Berker and Susan~R. McKay.
\newblock Hierarchical models and chaotic spin glasses.
\newblock {\em J. Statist. Phys.}, 36(5-6):787--793, 1984.

\bibitem{BleherLyubich1991}
P.~M. Bleher and M.~Yu. Lyubich.
\newblock Julia sets and complex singularities in hierarchical {I}sing models.
\newblock {\em Comm. Math. Phys.}, 141(3):453--474, 1991.

\bibitem{BleherZalis1989}
P.~M. Bleher and E~Zalis.
\newblock Asymptotics of the susceptibility for the {I}sing model on the
  hierarchical lattices.
\newblock {\em Commun. Math. Phys.}, 120:409--436, 1989.

\bibitem{BleherLyubichRoeder2017}
Pavel Bleher, Mikhail Lyubich, and Roland Roeder.
\newblock Lee-{Y}ang zeros for the {DHL} and 2{D} rational dynamics, {I}.
  {F}oliation of the physical cylinder.
\newblock {\em J. Math. Pures Appl. (9)}, 107(5):491--590, 2017.

\bibitem{BlokhBuffCheritatOversteegen2010}
A.~Blokh, X.~Buff, A.~Ch\'{e}ritat, and L.~Oversteegen.
\newblock The solar {J}ulia sets of basic quadratic {C}remer polynomials.
\newblock {\em Ergodic Theory Dynam. Systems}, 30(1):51--65, 2010.

\bibitem{BoscoRosa1987}
F.~A. Bosco and S.~{Gourlat Rosa jr.}
\newblock Fractal dimension of the {J}ulia set associated with the {Y}ang-{L}ee
  zeros of the {I}sing model on the {C}ayley tree.
\newblock {\em Europhys. Lett.}, 4(10):1103--1108, 1987.

\bibitem{Bowditch1998}
B.~H. Bowditch.
\newblock Markoff triples and quasi-{F}uchsian groups.
\newblock {\em Proc. London Math. Soc. (3)}, 77(3):697--736, 1998.

\bibitem{BrooksMatelski1978}
Robert Brooks and J.~Peter Matelski.
\newblock The dynamics of {$2$}-generator subgroups of {${\rm PSL}(2,\,{\bf
  C})$}.
\newblock In {\em Riemann surfaces and related topics: {P}roceedings of the
  1978 {S}tony {B}rook {C}onference ({S}tate {U}niv. {N}ew {Y}ork, {S}tony
  {B}rook, {N}.{Y}., 1978)}, volume~97 of {\em Ann. of Math. Stud.}, pages
  65--71. Princeton Univ. Press, Princeton, N.J., 1981.

\bibitem{BuffCheritat2012}
Xavier Buff and Arnaud Ch{\'e}ritat.
\newblock {Q}uadratic {J}ulia sets with positive area.
\newblock {\em Ann. Math. (2)}, 176(2):673--746, 2012.

\bibitem{BullettLomonaco2016ax}
S.~Bullett and L.~Lomonaco.
\newblock Mating quadratic maps with the modular group ii.
\newblock \url{https://arxiv.org/abs/1611.05257v1}, 2016.

\bibitem{BullettHaissinsky2007}
Shaun Bullett and Peter Ha\"{i}ssinsky.
\newblock Pinching holomorphic correspondences.
\newblock {\em Conform. Geom. Dyn.}, 11:65--89, 2007.

\bibitem{BullettPenrose1994}
Shaun Bullett and Christopher Penrose.
\newblock Mating quadratic maps with the modular group.
\newblock {\em Invent. Math.}, 115(3):483--511, 1994.

\bibitem{BullettPenrose1999}
Shaun Bullett and Christopher Penrose.
\newblock Perturbing circle-packing {K}leinian groups as correspondences.
\newblock {\em Nonlinearity}, 12(3):635--672, 1999.

\bibitem{Cao1995}
Chun Cao.
\newblock Some trace inequalities for discrete groups of {M}\"{o}bius
  transformations.
\newblock {\em Proc. Amer. Math. Soc.}, 123(12):3807--3815, 1995.

\bibitem{CarlesonGamelin1993}
Lennart Carleson and Theodore~W. Gamelin.
\newblock {\em Complex dynamics}.
\newblock Universitext: Tracts in Mathematics. Springer-Verlag, New York, 1993.

\bibitem{Cayley1879}
Professor Cayley.
\newblock Desiderata and {S}uggestions: {N}o. 3.{T}he {N}ewton-{F}ourier
  {I}maginary {P}roblem.
\newblock {\em Amer. J. Math.}, 2(1):97, 1879.

\bibitem{Crass2001}
Scott Crass.
\newblock Solving the quintic by iteration in three dimensions.
\newblock {\em Experiment. Math.}, 10(1):1--24, 2001.

\bibitem{Crass2002}
Scott Crass.
\newblock Solving the octic by iteration in six dimensions.
\newblock {\em Dyn. Syst.}, 17(2):151--186, 2002.

\bibitem{Crass2014}
Scott Crass.
\newblock Dynamics of a soccer ball.
\newblock {\em Exp. Math.}, 23(3):261--270, 2014.

\bibitem{DerridaDeSezeItzykson1983}
B.~Derrida, L.~de~Seze, C.~Itzykson, and and.
\newblock Fractal structure of zeros in hierarchical models.
\newblock {\em J. Statist. Phys.}, 33(3):559--569, 1983.

\bibitem{DerridaItzyksonLuck1984}
B.~Derrida, C.~Itzykson, and J.~M. Luck.
\newblock Oscillatory critical amplitudes in hierarchical models.
\newblock {\em Comm. Math. Phys.}, 94(1):115--132, 1984.

\bibitem{DouadyHubbard1985polylike}
Adrien Douady and John~Hamal Hubbard.
\newblock On the dynamics of polynomial-like mappings.
\newblock {\em Ann. Sci. \'Ecole Norm. Sup. (4)}, 18(2):287--343, 1985.

\bibitem{DoyleMcMullen1989}
Peter Doyle and Curt McMullen.
\newblock Solving the quintic by iteration.
\newblock {\em Acta Math.}, 163(3-4):151--180, 1989.

\bibitem{Eremenko2018talk}
Alexandre Eremenko.
\newblock Interactions between function theory and holomorphic dynamics.
\newblock \url{http://www.icms.org.uk/downloads/complex_talks/Eremenko.pdf},
  2018.
\newblock Resonances of complex dynamics.

\bibitem{GehringMartin1989}
F.~W. Gehring and G.~J. Martin.
\newblock Iteration theory and inequalities for {K}leinian groups.
\newblock {\em Bull. Amer. Math. Soc. (N.S.)}, 21(1):57--63, 1989.

\bibitem{GhulghazaryanAnanikyanJonassen2003}
Ruben Ghulghazaryan, Nerses Ananikyan, and Tore~M. Jonassen.
\newblock Julia sets and {Y}ang-{L}ee zeros of the {P}otts model on {B}ethe
  lattices.
\newblock In {\em Computational science---{ICCS} 2003. {P}art {I}}, volume 2657
  of {\em Lecture Notes in Comput. Sci.}, pages 85--94. Springer, Berlin, 2003.

\bibitem{Haruta1999}
Mako~E. Haruta.
\newblock Newton's method on the complex exponential function.
\newblock {\em Trans. Amer. Math. Soc.}, 351(6):2499--2513, 1999.

\bibitem{Haruta1992}
Mako~Emily Haruta.
\newblock {\em The dynamics of {N}ewton's method on the exponential function in
  the complex plane}.
\newblock ProQuest LLC, Ann Arbor, MI, 1992.
\newblock Thesis (Ph.D.)--Boston University.

\bibitem{HubbardSchleicherSutherland2001}
John Hubbard, Dierk Schleicher, and Scott Sutherland.
\newblock How to find all roots of complex polynomials by {N}ewton's method.
\newblock {\em Invent. Math.}, 146(1):1--33, 2001.

\bibitem{Hurley1986}
Mike Hurley.
\newblock Multiple attractors in {N}ewton's method.
\newblock {\em Ergodic Theory Dynam. Systems}, 6(4):561--569, 1986.

\bibitem{Jorgensen1976}
Troels J{\o}rgensen.
\newblock On discrete groups of {M}\"{o}bius transformations.
\newblock {\em Amer. J. Math.}, 98(3):739--749, 1976.

\bibitem{Jorgensen1977}
Troels J{\o}rgensen.
\newblock A note on subgroups of {$SL(2,{\bf C})$}.
\newblock {\em Quart. J. Math. Oxford Ser. (2)}, 28(110):209--211, 1977.

\bibitem{Keen1983}
Linda Keen.
\newblock Teichmueller spaces of punctured tori. {I}, {II}.
\newblock {\em Complex Variables Theory Appl.}, 2(2):199--211, 213--225, 1983.

\bibitem{KeenSeries1994}
Linda Keen and Caroline Series.
\newblock The {R}iley slice of {S}chottky space.
\newblock {\em Proc. London Math. Soc. (3)}, 69(1):72--90, 1994.

\bibitem{KhavinsonNeumann2006}
Dmitry Khavinson and Genevra Neumann.
\newblock On the number of zeros of certain rational harmonic functions.
\newblock {\em Proc. Amer. Math. Soc.}, 134(4):1077--1085, 2006.

\bibitem{KhavinsonNeumann2008}
Dmitry Khavinson and Genevra Neumann.
\newblock From the fundamental theorem of algebra to astrophysics: a
  ``harmonious'' path.
\newblock {\em Notices Amer. Math. Soc.}, 55(6):666--675, 2008.

\bibitem{KhavinsonSwiatek2003}
Dmitry Khavinson and Grzegorz \'{S}wi\c{a}tek.
\newblock On the number of zeros of certain harmonic polynomials.
\newblock {\em Proc. Amer. Math. Soc.}, 131(2):409--414, 2003.

\bibitem{KlimenkoKopteva2002}
E.~Klimenko and N.~Kopteva.
\newblock Discreteness criteria for {$\mathcal{R}\mathcal{P}$} groups.
\newblock {\em Israel J. Math.}, 128:247--265, 2002.

\bibitem{Klimenko2001}
Elena Klimenko.
\newblock Some examples of discrete groups and hyperbolic orbifolds of infinite
  volume.
\newblock {\em J. Lie Theory}, 11(2):491--503, 2001.

\bibitem{Lattes1918}
S.~{Latt\`es}.
\newblock {Sur l'it\'eration des substitutions rationnelles et les fonctions de
  {\it Poincar\'e}.}
\newblock {\em {C. R. Acad. Sci., Paris}}, 166:26--28, 1918.

\bibitem{ManeSadSullivan1983}
R.~Ma\~{n}\'{e}, P.~Sad, and D.~Sullivan.
\newblock On the dynamics of rational maps.
\newblock {\em Ann. Sci. \'{E}cole Norm. Sup. (4)}, 16(2):193--217, 1983.

\bibitem{Manning1992}
Anthony Manning.
\newblock How to be sure of finding a root of a complex polynomial using
  {N}ewton's method.
\newblock {\em Bol. Soc. Brasil. Mat. (N.S.)}, 22(2):157--177, 1992.

\bibitem{McKayBerkerKirkpatrick1982}
Susan~R. McKay, A.~Nihat Berker, and Scott Kirkpatrick.
\newblock Spin-glass behavior in frustrated {I}sing models with chaotic
  renormalization-group trajectories.
\newblock {\em Phys. Rev. Lett.}, 48(11):767--770, 1982.

\bibitem{McMullen1987}
Curt McMullen.
\newblock Families of rational maps and iterative root-finding algorithms.
\newblock {\em Ann. of Math. (2)}, 125(3):467--493, 1987.

\bibitem{McMullen1988}
Curt McMullen.
\newblock Braiding of the attractor and the failure of iterative algorithms.
\newblock {\em Invent. Math.}, 91(2):259--272, 1988.

\bibitem{McMullen1991}
Curt McMullen.
\newblock Rational maps and {K}leinian groups.
\newblock In {\em Proceedings of the {I}nternational {C}ongress of
  {M}athematicians, {V}ol. {I}, {II} ({K}yoto, 1990)}, pages 889--899. Math.
  Soc. Japan, Tokyo, 1991.

\bibitem{McMullen1994}
Curtis~T. McMullen.
\newblock The classification of conformal dynamical systems.
\newblock In {\em Current developments in mathematics, 1995 ({C}ambridge,
  {MA})}, pages 323--360. Int. Press, Cambridge, MA, 1994.

\bibitem{McMullen1994book}
Curtis~T. McMullen.
\newblock {\em Complex dynamics and renormalization}, volume 135 of {\em Annals
  of Mathematics Studies}.
\newblock Princeton University Press, Princeton, NJ, 1994.

\bibitem{Milnor2000}
John Milnor.
\newblock Periodic orbits, externals rays and the {M}andelbrot set: an
  expository account.
\newblock {\em Ast\'{e}risque}, (261):xiii, 277--333, 2000.
\newblock G\'{e}om\'{e}trie complexe et syst\`emes dynamiques (Orsay, 1995).

\bibitem{Milnor2006book}
John Milnor.
\newblock {\em Dynamics in one complex variable}, volume 160 of {\em Annals of
  Mathematics Studies}.
\newblock Princeton University Press, Princeton, NJ, third edition, 2006.

\bibitem{Przytycki2018ax}
F.~Przytycki.
\newblock Thermodynamic formalism methods in one-dimensional real and complex
  dynamics.
\newblock \url{https://arxiv.org/abs/1806.06186v1}, 2018.

\bibitem{Przytycki1989}
Feliks Przytycki.
\newblock Remarks on the simple connectedness of basins of sinks for iterations
  of rational maps.
\newblock In {\em Dynamical systems and ergodic theory ({W}arsaw, 1986)},
  volume~23 of {\em Banach Center Publ.}, pages 229--235. PWN, Warsaw, 1989.

\bibitem{RempeGillen2016ax}
Lasse Rempe-Gillen.
\newblock Arc-like continua, {J}ulia sets of entire functions and {E}remenko's
  conjecture.
\newblock \url{https://arxiv.org/abs/1610.06278v3}, 2018.

\bibitem{Rhie2001ax}
Sun~Hong Rhie.
\newblock Can a gravitational quadruple lens produce 17 images?
\newblock \url{https://arxiv.org/abs/astro-ph/0103463}, 2001.

\bibitem{Roeder2016}
Roland K.~W. Roeder.
\newblock Around the boundary of complex dynamics.
\newblock In {\em Dynamics done with your bare hands}, EMS Ser. Lect. Math.,
  pages 101--155. Eur. Math. Soc., Z\"{u}rich, 2016.

\bibitem{Schleicher2002}
Dierk Schleicher.
\newblock On the number of iterations of {N}ewton's method for complex
  polynomials.
\newblock {\em Ergodic Theory Dynam. Systems}, 22(3):935--945, 2002.

\bibitem{SchleicherStoll2017}
Dierk Schleicher and Robin Stoll.
\newblock Newton's method in practice: {F}inding all roots of polynomials of
  degree one million efficiently.
\newblock {\em Theoret. Comput. Sci.}, 681:146--166, 2017.

\bibitem{Shishikura2009}
Mitsuhiro Shishikura.
\newblock The connectivity of the {J}ulia set and fixed points.
\newblock In {\em Complex dynamics}, pages 257--276. A K Peters, Wellesley, MA,
  2009.

\bibitem{ShubSmale1986}
Michael Shub and Steve Smale.
\newblock On the existence of generally convergent algorithms.
\newblock {\em J. Complexity}, 2(1):2--11, 1986.

\bibitem{Smale1987}
Steve Smale.
\newblock On the efficiency of algorithms of analysis.
\newblock {\em Bull. Amer. Math. Soc. (N.S.)}, 13(2):87--121, 1985.

\bibitem{Sullivan1983}
Dennis Sullivan.
\newblock Conformal dynamical systems.
\newblock In {\em Geometric dynamics ({R}io de {J}aneiro, 1981)}, volume 1007
  of {\em Lecture Notes in Math.}, pages 725--752. Springer, Berlin, 1983.

\bibitem{Sullivan1985}
Dennis Sullivan.
\newblock Quasiconformal homeomorphisms and dynamics. {I}. {S}olution of the
  {F}atou-{J}ulia problem on wandering domains.
\newblock {\em Ann. of Math. (2)}, 122(3):401--418, 1985.

\bibitem{Tan1989}
Delin Tan.
\newblock On two-generator discrete groups of {M}\"{o}bius transformations.
\newblock {\em Proc. Amer. Math. Soc.}, 106(3):763--770, 1989.

\bibitem{Witt1990}
H.~J. Witt.
\newblock Investigation of high amplification events in light curves of
  gravitationally lensed quasars.
\newblock {\em Astron. Astrophys.}, 236:311--322, 1990.

\bibitem{Yoccoz2002}
Jean-Christophe Yoccoz.
\newblock Analytic linearization of circle diffeomorphisms.
\newblock In {\em Dynamical systems and small divisors ({C}etraro, 1998)},
  volume 1784 of {\em Lecture Notes in Math.}, pages 125--173. Springer,
  Berlin, 2002.

\end{thebibliography}
